\documentclass{article}

\usepackage[english]{babel}
\usepackage{mathtools}
\usepackage{relsize}

\usepackage[
backend=biber,
style=alphabetic,
sorting=ynt
]{biblatex}

\usepackage[letterpaper,top=3cm,bottom=3cm,left=3.5cm,right=3.5cm,marginparwidth=0cm]{geometry}

\usepackage{amsmath}
\usepackage{graphicx}
\usepackage[colorlinks=true, allcolors=blue]{hyperref}

\title{Generating Functions for Restricted Motzkin Paths in a Slit for Arbitrary Weights }
\author{
Cetin Hakimoglu-Brown  
\\
mathemails@proton.me
}
\begin{document}
\maketitle

\begin{abstract}
In this paper, I derive a generating function for discrete Motzkin paths of step size +2,+1,-1, restricted between two absorbing parallel lines, and generalized to arbitrary step weights. This continues off work of similar problems involving directed paths with either single barriers, uniform step size, or uniform weights.

\end{abstract}

\section{Introduction}

Motzkin paths and other variants of directed paths are well-studied objects in combinatorics. Most papers concern Motzkin paths with typically two, but not all, of the following properties: a single barrier, uniform step size, or uniform weights. Examples include   Feller (1957), Lengyel, T. (2009), Krattenthaler (2000, 2016), l-Shehawey (2008), and  Oste and Van der Jeugt (2015). Path enumerating generating functions involving double barriers (e.g. a slit) and uneven steps and varying weights are harder to come by, such as by C. Banderier (2010) and  Khalid and Prellberg (2019). The formulas are considerably more complicated, and analytic expressions tend to have a considerable number of terms, or are un-simplified. This paper simplifies a expression for a +1, +2, -1 Motzkin path between a slit and arbitrary weights, to derive an elegant expression similar to well-known expressions involving Chebyshev polynomials, but with double summations. 

\subsection{Problem Statement}

Khalid and Prellberg [11] give an analytic expression for a double-barrier +2,+1,-1 step-size  Motzkin Path for arbitrary weights, but it's very cumbersome to write out, involving the three roots of a cubic of characteristic equation:

\begin{multline}
G_{(u,v)}^{w,2,1}(t)=
-\frac{1}{tp_2}\times \\
\dfrac{\mathlarger{\mathlarger{\mathlarger{\sum}}}\limits_{l=0}^r\left(
   \splitdfrac{\splitdfrac{z_1^{w+2}(z_2^{w-(v-u)_+-l+1}z_3^{(u-v)_++l}-z_3^{w-(v-u)_+-l+1}z_2^{(u-v)_++l})}{- z_2^{w+2}(z_1^{w-(v-u)_+-l+1}z_3^{(u-v)_++l}-z_3^{w-(v-u)_+-l+1}z_1^{(u-v)_++l})}}{+z_3^{w+2}(z_1^{w-(v-u)_+-l+1}z_2^{(u-v)_++l}-z_2^{w-(v-u)_+-l+1}z_1^{(u-v)_++l})}
\right)}
{z_1^{w+3}(z_2^{w+2}-z_3^{w+2})-z_2^{w+3}(z_1^{w+2}-z_3^{w+2})+z_3^{w+3}(z_1^{w+2}-z_2^{w+2})}
\end{multline}

Using a toolbox of  various analytical combinatorics techniques, I will find a simplification that produces a nice closed-form expression for non-trivial initial conditions. 

\section{Derivation of Generating Function}\label{sec:2} 

\subsection{Characteristic Equation}

Consider that a particle has three degrees/weights, which can be probabilities. The characteristic equation for the three-term recursion for a +2,+1,-1 walk is is given by:

\begin{equation}
\frac{a_3x^3}{a_1}+\frac{a_2x^2}{a_1}-\frac{x}{za_1}+q=0
\end{equation}

In which weights $a_3,a_2,a_1$ are defined:

$a_3$ = forward 2 steps

$a_2$ = forward 1 step

$a_1$ = backward 1 step

This cubic is quite general, and the final generating function will be expressed in terms of $z$, but the $q$ term is used later in the derivation, but is not part of the generating function.

This is covered in many elementary/intro texts, such as [12].

\subsection{Boundary Conditions and System of Equations}

Using the procedure in [11], without Schur functions, the generating function arises from a system of three linear equations that encode the boundary conditions. The simplest non-trivial case involves the particle starting at $s=m-1$ and we seek to find the entry corresponding to powers of the associated matrix. This can either be 'path steps' of probabilities of absorption on the left-sided barrier conditional on having not been absorbed before.  

The process also is described by the associated $m+2$ by $m+2$ heptadiagonal transition matrix (in this example, $m=5$). Here we can see absorption occurs if the particle hits the squares corresponding to $0$, $m$, or $m+1$: 

\begin{displaymath}
\begin{bmatrix}
1 & 0 & 0 & 0 & 0 & 0 & 0 \\
a_1 & 0 & a_2 & a_3 & 0 & 0 & 0 \\
0 & a_1 & 0 & a_2 & a_3  & 0 & 0 \\
0 & 0 &  a_1 & 0 & a_2 & a_3  & 0 \\
0 & 0 & 0 & a_1  & 0 & a_2 & a_3  \\
0 & 0 & 0 & 0 & 0 & 1 & 0 \\
0 & 0 & 0 & 0 & 0 & 0 & 1 
\end{bmatrix}
\end{displaymath}

The system of linear equations and the resulting generating function, in which $a,b,c$ are the roots of the cubic (2):

\begin{equation}
    \begin{cases}
      x+y+z=1\\
      xa^{m}+yb^{m}+zc^{m}=0\\
      xa^{m+1}+yb^{m+1}+zc^{m+1}=0
    \end{cases}\,.
\end{equation}

\begin{equation}
    \begin{cases} 
    
    x =  \frac{(b-c) b^m c^m}{(a^{m + 1} (b^m - c^m) + a^m (c^{m + 1} - b^{m + 1}) + (b - c) b^m c^m)}\\
    y =  \frac{(a-c) a^m c^m}{-(a^{m + 1} (b^m - c^m) - a^m (c^{m + 1} - b^{m + 1}) - (b - c) b^m c^m)}\\
    z =  \frac{(a-b) b^m a^m}{(a^{m + 1} (b^m - c^m) + a^m (c^{m + 1} - b^{m + 1}) + (b - c) b^m c^m)}\\ 
    
     \end{cases}\,.
\end{equation}

The probability generating function is thus of the form: 

\begin{equation}
f(z,s,q,m)=xa^s+yb^s+zc^s ; m+1 \geq s\geq 0
\end{equation}

Which is similar to  (1).

The exact generating function for the simplest but non-trivial case setting $s=m-1$ corresponding to (5) for the exact probability of absorption on the left-side barrier, conditional on having not been absorbed before, simplifies to:

\begin{equation}
f(z,q,m) = \frac{-a_3}{qa_1}\left(\frac{(a-b)(b-c)(c-a)}{(b-c)a^{-m}+(a-b)c^{-m}+(c-a)b^{-m}} \right) 
\end{equation}
 
\subsection{Lagrange Inversion Theorem}

We generate series approximation for roots $a,b,c$ of the cubic (2). Let $b$ be the small root and let $a,c$ be the large ones. 

The small root behaves as $b\approx a_1zq$. The large roots behave as $a\approx R_1,c\approx R_2$

Because we are working with formal power series, ensuring convergence is not necessary. 

We perform the following steps:

1. For root $b$, we apply the Lagrange inversion theorem to the auxiliary quadratic, which has two roots $R_1, R_2$  in terms of $a_2, a_3$:

\begin{equation}
\frac{a_3x^2}{a_1}+\frac{a_2x}{a_1}-\frac{1}{za_1}=0 
\end{equation}

A generalized version of the Lagrange inversion theorem is applied that allows for powers $m$ of roots to be expressed as a single series.

2.The Cauchy product is applied to step 1 to convert the product of two summations to a single summation containing a product of two binomials, allowing for terms to be extracted more easily in a closed-form expression. This is discussed in more detail in [5]. 

3. We observe that the resulting polynomial of step 2, in terms of $z$, has symmetric properties, allowing the number of terms summed to be cut in half, denoted by the floor function.
\newpage
Applying steps 1-3, define:

\begin{multline}
g_n=\cos\left( \frac{n\pi}{2}\right)  \binom{m-1-n}{n/2} ^{2}+ \\ \sum_{u=0}^{\left \lfloor \frac{n-1}{2} \right \rfloor} (-1)^{n+u}\binom{m-1-n}{u} \binom{m-1-n}{n-u}  (R_1^{2u-n}+R_2^{2u-n} )  
\end{multline}

Now define: 
\begin{equation}
R_1=g+hi, R_2=g-hi ,g=\frac{-a_2}{2a_3}, h=\frac{(a_2^{2}+4a_3/z)^{1/2}}{2a_3}
\end{equation}

We observe that for positive weights, that the argument lies in the second quadrant, so after performing some labor we  have:

\begin{multline}
g_n=\cos\left( \frac{n\pi}{2}\right)  \binom{m-1-n}{n/2} ^{2}+\\  2\sum_{u=0}^{\left \lfloor \frac{n-1}{2} \right \rfloor} (-1)^{5n/2}\binom{m-1-n}{u} \binom{m-1-n}{n-u} \cos \left( (2u-n)\arccos {\left( \frac{ia_2\sqrt{z}}{2\sqrt{a_3}} \right)}\right) 
\end{multline}

Putting it all together, we have the series solution for one of the roots of the cubic:

\begin{equation}
b^{1-m}=(a_1zq)^{1-m}  \left( 1+(1-m)  \sum_{n=1}^{}{\frac{(a_1^2a_3)^{n/2}z^{3n/2}q^ng_n}{n-m+1}}             \right) 
\end{equation}

By using the multiple-angle formula for cosines, which is related to the Chebyshev polynomials, we can extract powers of $z$. The complex and fractional powers of $z$ and $a_3$ cancel out in such a way as to ensure that the above series is rational and real, as expected for a formal power series.

And applying step 1, obtain the approximations for the other two roots (the full series not being needed for the proof):

\begin{equation}
a^{1-m}\approx R_1^{-m+1}+\left( \frac{(1-m)R_1^{-m-1}q}{R_1-R_2}  \right) 
\end{equation}

\begin{equation}
c^{1-m}\approx R_2^{-m+1}+\left( \frac{(1-m)R_2^{-m-1}q}{R_2-R_1}  \right) 
\end{equation}
\newpage
\subsection{Reduction Theorem}

In this section I will introduce a theorem that allows for the simplification of Schur and Vandermonde-like polynomials that arises from computing determinants. I have yet to see this method employed in quite the same way for enumeration problems.  For bounded directed paths  beyond just symmetric step-size (+1,-1), the generating function  takes on a more complicated form in which the numerator and denominator contain products of the roots, such as for a cubic, and this method helps simply such products.

Let $F(x,q) = b_0x^n+b_1x^{n-1}+...+q$, such as the characteristic equation given by (2).

Then for each root $r_i$ for $i=1,2..n$, the identities below hold. For example for $r_1$:

\begin{equation}
  \epsilon \frac{r_1^m}{b_0(r_1-r_2)(r_1-r_3)(r_1-r_4)...}=\frac{r_1^{m+1}}{m+1}\frac{d}{dq}
\end{equation}

also

\begin{equation}
 \epsilon(r_1^m) \prod_{2 \leq j < k \leq n} (r_j-r_k)=b_0^{2-n}\sqrt{D}\left[\frac{r_1^{m+1}}{m+1}\frac{d}{dq}\right]
\end{equation}

In which $\epsilon =  \pm 1$, $D$ is the discriminant of $F(x,q)$,  and $r^{m}_i$  is the series expansion of the specified root about $q$.

Proof:

(14) and (16) are the same, so we'll work on (14). 

Step 1. We have: $\frac{r_1^{m+1}}{m+1}\frac{d}{dq}=r_1^{m} \left[r_1\frac{d}{dq}\right] $.

Step 2. Let $  F(x,q)=b_0(x-r_1)(x-r_2)...$  Via Vieta's formulas 

Step 3. Let ${x \to r_1}$, take the reciprocal of both sides and we have:

\begin{equation}
\lim_{x \to r_1} \frac{F(x,q)}{x-r_1}= \left[ r_1\frac{d}{dq}\right]^{-1}
\end{equation}

Step 4. Computing the Laurent series of $\frac{F(x,q)}{x-r_1}$ about $x=r_1$ and extracting the residue, we have: 

\begin{equation}
b_0(r_1-r_2)(r_1-r_3)...=\left.F(x,q)\frac{d}{dx}\right|_{F^{-1}(q)=r_1} =\left[F^{-1}(q)\frac{d}{dq}\right]^{-1}
\end{equation}

This is the well-known formula for the derivative of an inverse function, completing the proof.
\newpage
\subsection{Simplification}

Using (15) and the fact $ b_0^{-2} \sqrt{D}=-(a-b)(b-c)(c-a)$ (choosing the minus sign), letting $b_0=a_3/a_1$ and $m=-m$ and $q=1$, we can simplify (6) as:

\begin{equation}
f(z,m) =  \left[ \left[\frac{b^{1-m}}{1-m} +\frac{a^{1-m}}{1-m}+\left.\frac{c^{1-m}}{1-m} \right]\frac{d}{dq} \right|_{q=1} \right]^{-1} 
\end{equation}

Plugging in (11),(12),and (13) into (18) we have: 

\begin{equation}
f(z,m) = z^{m} \left[  za_1^{1-m} \left[1  + \sum_{n=1}^{}{(a_1^2a_3)^{n/2}z^{3n/2}g_n}      \right] +z^{m} \left[  \frac{R_1^{-1-m}-R_2^{-1-m}}{R_1-R_2} \right] \right]^{-1} 
\end{equation}

By the Faddeev–LeVerrier algorithm, the denominator of (19) is a polynomial in terms of $z$ is of degree $m-1$, and the corresponding transition matrix has a maximum rank of $m-1$ (omitting the trivial eigenvalues 0,1). 
 
Expanding the Puiseux series  of  $\frac{R_1^{-1-m}-R_2^{-1-m}}{R_1-R_2}$  about z, we find that the lowest  power of z grows by $z^{(m+2)/2}$  for m even and $z^{(m+3)/2}$  for m odd, which when  multiplied by $z^{m}$  is obviously of a higher power than $z^{m-1}$, and hence does not contribute to the generating function   and   we can discard  it.

Finally we have: 

\begin{equation}
f(z,m) = z^{m-1}a_1^{m-1}   \left[1  + \sum_{n=1}^{}{(a_1^2a_3)^{n/2}z^{3n/2}g_n}      \right]^{-1} 
\end{equation}

This completes the proof. What we have done is simplified an expression involving multiple root extractions, powers of roots, and products of roots of a generalized cubic into a single summation that involves a second summation. 

\subsection{Example}

Let $a_3=2,a_2=3,a_1=1,m=9$. Computing the summation we have:

\begin{align}
n=1 &= -21z^2 \\
n=2 &= -12z^3+135z^4 \\
n=3 &= 120z^5-270z^6 \\
n=4 &= 81z^8-216z^7 + 24z^6\\
n=5 &= -36z^8 
\end{align}

The generating function is thus:
\begin{equation}
f(z,9)=\frac{z^8}{1 - 21 z^2 - 12 z^3 + 135 z^4 + 120 z^5 - 246 z^6 - 216 x^7 + 45 x^8}
\end{equation}

Which has the expansion $\sum_{n=m-1}^{k} c_n z^n$: 

\begin{multline}
z^8 + 21 z^{10} + 12 z^{11} + 306 z^{12} + 384 z^{13} + 3981 z^{14} + 7812 z^{15}\\ + 50580 z^{16} + 130752 z^{17} + 649332 z^{18} + 1980432 z^{19} + 8487756 z^{20}...
\end{multline}

In general for (20) and assuming $a_1=1$, the expansion $c_n$ is also is equal to the $A_{(m,2)}$ entry of powers of the associated $m+2$ by $m+2$ sized absorbing transition matrix, such that:

\begin{equation}
f(z,m) = \sum_{n=1}^{k} c_n z^n,  = \sum_{n=1}^{k} A_{(m,2)}^{n-1}  z^{n}
\end{equation}

\subsection{Omitting the +1 Step}

If $a_2=0$, then the generating function simplifies considerably. The additional degree of freedom of allowing a +1 step introduces vastly more complexity into the final expression. For (10), extracting powers of $z$ from the cosine component requires yet another summation, for a total of three summations compared to a single summation if $a_2=0$. The  generating  function is thus, assuming a starting position $s=m-1$ and the same boundary conditions as earlier: 

\begin{equation}
f(z,m) = z^{m-1}a_1^{m-1}   \left[1  + \sum_{n=1}^{}{(a_1^2a_3)^{n}z^{3n}\binom{3n-m}{n}}     \right]^{-1};3n-m<0 
\end{equation}

Via the Vandermonde's identity or by using the Lagrange Inversion Theorem again for a trinomial cubic, letting $a_3=a_1=1$ and equating the coefficients $z$, the binomials given by $z^{3n}g_{2n}$ and $z^{3n}\binom{3n-m}{n}$ are equal.

\subsection{Summation Bounds on Generating Function}

For formula (20), it's reasonable to inquire as to how many values of $n$ of $g_n$ must be calculated to compute the necessary terms of the generating function, up to $z^{m-1}$. For example, 2.6 for $m=9$, we have five terms. I conjecture the exact minimum is given by:

\begin{equation}
{\left \lfloor \frac{2(m-1)}{3} \right \rfloor}
\end{equation}

An additional $g_{n}$ term produces a leading  exponent of $z^m$ (with a zero-valued coefficient), which is extraneous given that coefficients only up to $z^{m-1}$ are needed.

To prove this, let $m=3v+1,m=3v+2,m=3v$ for some positive integer $v>0$

We can express the three possible cases as a table:
\begin{table}[h!]
\centering
 \begin{tabular}{||c c c c||} 
 \hline
Case & Sum1 & n (max) & m \\ [0.5ex] 
 \hline\hline
1 & v-1 & 2v & 3v+2 \\ 
 2 & v-1 & 2v & 3v+1 \\
 3 & v-1 & 2v-1 & 3v \\
  [1ex] 
 \hline
 \end{tabular}
\end{table}

These values are derived by plugging values of column $m$ into (30) to get the maximum $n_{max}$ . The Sum1 column is from plugging values of column $n_{max}$ into $\left \lfloor \frac{n-1}{2} \right \rfloor$, which is given by (10).
\newpage
For case 1, $m=3v+2,n=2v$, so we have:

\begin{multline}
g_nz^{3v}=z^{3v}\cos\left( v\pi\right)  \binom{v+1}{v} ^{2}+\\  z^{3v}2\sum_{u=0}^{v-1} (-1)^{5v}\binom{v+1}{u} \binom{v+1}{2v-u} \cos \left( (2u-2v)\arccos {\left( \frac{ia_2\sqrt{z}}{2\sqrt{a_3}} \right)}\right) 
\end{multline}

$ \binom{v+1}{2v-u}=0 $ for $0 \le u  < v-1$ and 1 when  $u  = v-1$

Letting $u=v-1$ , we have $\cos  (2\arccos {(ik z^{1/2})}) = -2zk^2-1$ for some constant $k$. Hence we see we have a $z^{3v+1}$ leading exponent, which is $z^{m-1}$.

For case 2, $m=3v+1,n=2v$, so we have:

\begin{multline}
g_nz^{3v}=z^{3v}\cos\left( v\pi\right)  \binom{v}{v} ^{2}+\\  z^{3v}2\sum_{u=0}^{v-1} (-1)^{5v}\binom{v}{u} \binom{v}{2v-u} \cos \left( (2u-2v)\arccos {\left( \frac{ia_2\sqrt{z}}{2\sqrt{a_3}} \right)}\right) 
\end{multline}

Here we see  $ \binom{v}{2v-u}=0 $ for $0 \le u  \le v-1$, and that the leading exponent is  $z^{3v}$, which is $z^{m-1}$.

For case 3, $m=3v,n=2v-1$, so we have:

\begin{multline}
g_nz^{3v-3/2}=0+  z^{3v-3/2}2\sum_{u=0}^{v-1} (-1)^{5v-5/2}\binom{v}{u} \binom{v}{2v-1-u} \cos \left( (2u-2v+1)\arccos {\left( \frac{ia_2\sqrt{z}}{2\sqrt{a_3}} \right)}\right) 
\end{multline}

$\cos\left((\pi/2)(2v-1)\right)=0$. Here we see  $ \binom{v}{2v-1-u}=0 $ for $0 \le u  < v-1$,  and 1 when  $u  = v-1$. And $\cos ((2u-2v+1) \arccos (ikz^{1/2})) = ikz^{1/2}$ for constant $k$ when $u =  v-1$.  Also $ikz^{1/2} (-1)^{5v-5/2}$ is real. Similarly, the fractional exponent of $a_3$ becomes an integer value.  Thus, the leading exponent is  $z^{3v-1}$, which is $z^{m-1}$.

 

\newpage

\begin{thebibliography}{9}
\bibitem{}
Helmut Prodinger, On k-Dyck paths with a negative boundary, 2020. https://arxiv.org/pdf/1912.06930.pdf

\bibitem{}
WILLIAM FELLER (1957) An Introduction
to Probability Theory
and Its Applications
, John Wiley and Sons, Inc.

\bibitem{}
Christian Krattenthaler (2017) Lattice Path Enumeration, 
, https://arxiv.org/pdf/1503.05930.pdf


\bibitem{}
Cyril Banderier, Philippe Flajolet  (2002) Basic analytic combinatorics of directed lattice paths, https://www.sciencedirect.com/science/article/pii/S0304397502000075

\bibitem{}
Cyril Banderier, Christian Krattenthaler, Alan Krinik, Dmitry Kruchinin, Vladimir Kruchinin, et al..
Explicit formulas for enumeration of lattice paths: basketball and the kernel method. Developments in
Mathematics, Springer, 2019, Lattice Path Combinatorics and Applications, pp.78-118. ffhal-01456197f

\bibitem{}
Tamas Lengyel (2009), The conditional gambler’s ruin problem with ties allowed, https://core.ac.uk/download/pdf/82111895.pdf

\bibitem{}
 Christian Krattenthaler (2000) Permutations with restricted patterns and Dyck paths, https://arxiv.org/abs/math/0002200
 
\bibitem{}
M. A. El-Shehawey (2008) Trinomial random walk with one or two imperfect absorbing barriers https://link.springer.com/article/10.2478/s12175-008-0080-5


\bibitem{}
Anum Khalid and Thomas Prellberg
 (2020) Skew Schur Function Representation of Directed Paths in a Slit,
 https://arxiv.org/abs/1907.09842

\bibitem{}
Roy Oste and Joris Van der Jeugt (2015) Motzkin paths, Motzkin polynomials and recurrence relations
 
\bibitem{}
Cetin Hakimoglu (2022) Generating Functions for Asymmetric Random Walk Processes With Double Absorbing Barriers, https://arxiv.org/abs/2212.09919 


\bibitem{}

Marie-Louise Lackner Michael Wallner (2015) An invitation to analytic combinatorics and lattice path counting, https://dmg.tuwien.ac.at/mwallner/files/lpintro.pdf



\end{thebibliography}
\end{document}